\documentclass[letter,9pt]{beamer}

\usepackage{ifthen}
\usepackage[T1]{fontenc}
\usepackage[ansinew]{inputenc}
\usepackage{geometry}
\usepackage{verbatim} 
\usepackage{color}
\usepackage{ragged2e}
\usepackage{newlfont}
\usepackage{pifont}

\usepackage{yfonts}
\usepackage{dsfont}
\usepackage{stmaryrd} 
\usepackage{bbm} 
\usepackage{xspace}
\usepackage{amsmath,amssymb,amsfonts,textcomp}
\RequirePackage{amsthm}
\usepackage[mathscr]{eucal}
\usepackage{yhmath}
\RequirePackage{multirow}
\usepackage{pgf} %
\usepgfmodule[matrix]
\usepackage{tikz}

\usetikzlibrary{%
  calc,%
  arrows,%
  trees,%
  shapes.misc,
  shapes.arrows,%
  shapes.geometric,%
  chains,%
  matrix,%
  positioning,
  scopes,%
decorations.markings,%
  decorations.pathreplacing,%
  decorations.pathmorphing,
  shadows,%
backgrounds%
}

\tikzset{
grph/.style ={>=latex', shorten <=-.3em, shorten >=-.3em},%
morph/.style ={>=stealth',shorten <=-.3em, shorten >=-.3em},%
every on chain/.append style={join},%
every join/.style={->},%
/seqStyle/.style ={>=stealth'},%
}

\pgfkeys{
/Arrowstyle/.code=\pgfkeys{/Arrowstyle/.initial=#1},/Arrowstyle=morph,
/Dir/.code=\pgfkeys{/Dir/.initial=#1},/Dir=->,
/Dist/.code=\pgfkeys{/Dist/.initial=#1},/Dist=2,
/Bline/.code=\pgfkeys{/Bline/.initial=#1}, /Bline=0.0em,
/Lab_level/.code=\pgfkeys{/Lab_level/.initial=#1},/Lab_level=0,
/arrowstyle/.initial=morph,/dir/.initial=->,/above/.initial={},
/below/.initial={},/lab_level/.initial=0,/dist/.initial=2,/bline/.initial=-0.0em,
/isomark/.initial={\sim},/mark_level/.initial=-0.12,/side/.initial=above,
/style/.initial=, /syst/.initial={}%
}

\newtheorem{Thm}{Theorem}[section]
\newtheorem{Cor}[Thm]{Corollary}
\newtheorem{Lem}[Thm]{Lemma}
\newtheorem{Prop}[Thm]{Proposition}

\theoremstyle{definition}
\newtheorem{Defn}{Definition}[section]
\newtheorem{Rem}[Defn]{Remark}

\newenvironment{prv}{
\begin{proof}}{\end{proof} } 

\newcommand{\ParIt}[1]{ \mbox{ } \\ \smallskip \noindent\textit{#1}.~}

\numberwithin{equation}{section}

\newcommand*{\midsloppy}{%
  \tolerance 5000%
  \hbadness 4000%
  \emergencystretch 1.5em%
  \hfuzz .1  
  \vfuzz\hfuzz}

\newenvironment{cent}[1][\vskip-0.8em]{ #1
\begin{center}}{ \end{center} 
\vskip-0.2em }

\newcommand{\myeqar}[2][rl]{\begin{array}{#1} #2 \end{array}}

\newcommand*{\nsc}[1][\hspace{-0.04em}]{#1,#1}

\newcommand*{\mrm}[1]{\mathrm{#1}}
\newcommand*{\msf}[1]{\mathsf{#1}}
\newcommand*{\mfr}[1]{\mathfrak{#1}}

\newcommand*{\sscr}[1]{{\scriptscriptstyle{#1}}} 

\newcommand*{\bbar}[1]{\overline{#1}} 
\newcommand*{\undl}[1]{\underline{#1}}   
\newcommand*{\wtilde}[1]{\widetilde{#1}} 
\newcommand*{\prim}{^\sscr{\prime}} 

\newcommand*{\inj}[1][\hspace{0.1em}]{
\inds{i\hspace{-0.1em}}{#1}} 
\newcommand*{\proj}[1][\hspace{0.1em}]{
\inds{p\hspace{-0.1em}}{#1}} 
\newcommand*{\iinj}{\mrm{i}\hspace{-0.2em}} 
 
\newcommand*{\phiso}[1][]{\phi_{{\scriptscriptstyle{#1}}}}

\newcommand*{\K}{\mathsf{\mrm{K}}}

\newcommand*{\kk}{\mathsf{\mrm{k}}}

\newcommand*{\C}{\mathcal{C}}

\newcommand*{\T}{\cal{T}}

\DeclareMathOperator{\add}{\mrm{add}}
\DeclareMathOperator{\Ind}{\mrm{ind}}

\DeclareMathOperator{\Dual}{\mathrm{D}}

\newcommand*{\Hm}[1][]{\mrm{Hom}_{{#1}}}

\newcommand*{\Hc}[1]{\ind{\mrm{Hom}}{\C}(#1)}

\newcommand*{\HM}[2]{\ind{\mrm{Hom}}{#1}(#2)}

\newcommand*{\J}[2][]{\indu{\mrm{J}}{#2}{#1}}

\newcommand*{\Irr}[1][]{\mrm{Irr}_{#1}}

\newcommand*{\Ext}[2][1]{\indu{\mrm{Ext}}{#2}{#1}}

\newcommand*{\B}{\mathrm{B}}

\newcommand{\dm}[1][n]{\mrm{d}}

\newcommand*{\set}[1]{\left\{#1\right\}}

\newcommand*{\Psmatr}[1]{\left[ \hspace{-0.1em}\begin{smallmatrix}
#1 \end{smallmatrix} \hspace{-0.1em}\right]} 

\newcommand*{\ind}[2]{{#1\hspace{-0.02em}}_{#2}}
\newcommand*{\inddl}[3]{{_{#2\hspace{-0.02em}}}{#1\hspace{-0.02em}}_{#3}}
\newcommand*{\inds}[2]{{#1\hspace{-0.02em}}_{{\scriptscriptstyle{#2}}}}
\newcommand*{\indu}[3]{{#1\hspace{-0.02em}}_{#2}^{#3}}

\newcommand*{\Som}[2][]{\sum\limits_{#2}^{#1}\hspace{-0.12ex}}

\newcommand*{\sdash}{\textrm{-}}

\newcommand*{\ssminus}{{\scriptstyle{\smallsetminus}}}

\newcommand*{\vphi}[1][{\hspace{0.015em}}]{\ind{\varphi}{#1}}

\newcommand*{\tr}[1][]{\mfr{t}_{#1}}

\newcommand*{\bilbox}[1][]{\Box}

\newcommand*{\bss}[1][]{\inds{b}{#1\hspace{0.02em}}}

\newcommand*{\Bimd}[3][]{{\inddl{B}{#2}{#3}^{#1}}}

\newcommand*{\op}[1]{#1^{{\scriptscriptstyle{\circ}}}}

\renewcommand*{\to}[1][]{\hspace{-0.4em}
\pgfkeys{/arrowstyle=\pgfkeysvalueof{/Arrowstyle},/dir=\pgfkeysvalueof{/Dir},
/style=,/dist=\pgfkeysvalueof{/Dist},/lab_level=\pgfkeysvalueof{/Lab_level},/above={},
/below={},/mark_level=-0.12,/bline=-0.25em,/isomark={},/side=above,#1}
\begin{tikzpicture}[\pgfkeysvalueof{/arrowstyle},baseline=\pgfkeysvalueof{/bline},
node distance=\pgfkeysvalueof{/dist}em]
\node (i) {}; \node (j)[right=of i]  {}; \node (label_pos) [inner
sep=0.12em] at ($ (i.east)!0.5! (j.west) $) {};
\path[\pgfkeysvalueof{/dir},\pgfkeysvalueof{/style}] (i) edge (j);
\node [inner sep=0em,above=\pgfkeysvalueof{/lab_level}em of
label_pos] {{$\scriptstyle{\pgfkeysvalueof{/above}}$}};
\node [inner sep=0em,below=\pgfkeysvalueof{/lab_level}em of
label_pos] {{$\scriptstyle{\pgfkeysvalueof{/below}}$}};
\node [inner sep=0em,
\pgfkeysvalueof{/side}=\pgfkeysvalueof{/mark_level}em of
label_pos] {${\scriptscriptstyle{\pgfkeysvalueof{/isomark}}}$};
\end{tikzpicture}
\pgfkeys{/above={},/below={},/mark_level=-0.12}\hspace{-0.4em}}

\newcommand*{\eqto}[1][\pgfkeysvalueof{/dist}]{
\pgftext{\begin{tikzpicture}[node distance=#1em] \node (pos1) {};
\node (pos2) [right=of pos1] {}; \path [shorten <=0.0em,shorten
>=0.0em,-] ($(0em,0.2em) + (pos1)$) edge ($(0em,0.2em) + (pos2)$)
($(0em,-0.2em) + (pos1)$) edge ($(0em,-0.2em) + (pos2)$);
\end{tikzpicture}}}

\newcommand*{\Seq}[2][]{
\pgfkeys{/Arrowstyle=morph,/Dist=1.5,/Lab_level=0.0,Dir=->,/bline=-0.2em,#1}
\makebox{\ensuremath{#2}}
\pgfkeys{/Arrowstyle=morph,/Dist=2,/Dir=->,/Lab_level=0}}
\newcommand*{\morph}[3][]{ 
\Seq{#2 \to[/dist=2,#1] #3}}
\newcommand*{\diagram}[3][]{
\pgfkeys{/Dist=2,/dist=1,clsep/.initial=2,/bline=-0.2em,#1}
\begin{tikzpicture}[morph,baseline=\pgfkeysvalueof{/bline}]
\matrix (m) [matrix of math nodes,ampersand replacement=\&, row
sep=\pgfkeysvalueof{/dist}em,column
sep=\pgfkeysvalueof{/clsep}em,text height=1em, text depth=0.4em]
{#2 \\}; #3;
\end{tikzpicture} \pgfkeys{/Dist=2}
}

\newcommand*{\trgl}[4][]{
\pgfkeys{/Arrowstyle=morph,/Dist=1.5,/Lab_level=0.0,Dir=->,/bline=-0.2em,
/u1/.initial={},/u2/.initial={},/u3/.initial={},/shiftobj/.initial={#2[1]},#1}
\makebox{\ensuremath{#2 \to[/above={\pgfkeysvalueof{/u1}}] #3
\to[/above={\pgfkeysvalueof{/u2}}] #4
\to[/above={\pgfkeysvalueof{/u3}}] \pgfkeysvalueof{/shiftobj}}}
\pgfkeys{/Arrowstyle=morph,/Dist=2,/Dir=->,/Lab_level=0,/u1/.initial={},
/u2/.initial={},/u3/.initial={}}}

\makeatletter

 \newcommand\makebeamertitle{\frame{\maketitle}}%
 \AtBeginDocument{%
   \let\origtableofcontents=\tableofcontents
   \def\tableofcontents{\@ifnextchar[{\origtableofcontents}{\gobbletableofcontents}}
   \def\gobbletableofcontents#1{\origtableofcontents}
 }

\allowdisplaybreaks[1]
\setbeamercovered{transparent}

\makeatother

\begin{document}

\title[Non-simply laced cluster structures]{A  non-simply laced version for cluster structures on $2$-Calabi-Yau categories}

\author[B. Nguefack]{Bertrand Nguefack\inst{1}}

\institute[Univ Yaounde I]{
\inst{1}University of Yaounde I; b.nguefack@uy1.uninet.cm \, 
ngc.bertrand@gmail.com
}

\date[December 2013]{
Beamer version of the  article published online at \href{
http://www.sciencedirect.com/science/article/pii/S0022404913002314}{\alert{Journal of Pure and Applied Algebra}} \, \, \href{http://dx.doi.org/10.1016/j.jpaa.2013.11.027}{\alert{DOI: 10.1016/j.jpaa.2013.11.027}} 
}

\makebeamertitle

\section{Extended abstract}
\begin{frame}{Sommary}

\tableofcontents{}

This paper investigates a non simply-laced version of cluster structures for  $2$-Calabi-Yau or stably $2$-Calabi-Yau categories over
arbitrary fields.  It results that   $2$-Calabi-Yau or stably 2-Calabi-Yau categories having a cluster tilting subcategory with neither loops nor 2-cycles do have the generalized version of cluster structure. This is in particular the case of cluster categories over non-algebraically closed fields.

\begin{block}{Keywords}
Calabi-Yau category, cluster structure, modulated quivers,  cluster algebra.

\textbf{MSC} (2010) Primary: 16G70;  Secondary: 18E30, 18E10,   13F60, 12E15 
\end{block}

\end{frame}

\begin{frame}{Brief introduction}

The theory of cluster algebras, introduced by Fomin and Zelevinsky,  is connected with many areas of mathematics. In the representation theory of algebras, the philosophy has been to provide a categorical interpretation of the main combinatorics, called mutation
of seeds, used to define cluster algebras. 

Most of the connections between cluster algebras and representation theory are proved only in the
simply-laced case which corresponds to cluster algebras associated with skew-symmetric matrices and to representations of quivers. The general case deals with skew-symmetrizable matrices and one should consider
categories over arbitrary fields, and hence work with representations of modulated quivers.

In this paper and some of our forthcoming work (\href{http://arxiv.org/abs/1004.2213}{\alert{Potentials and Jacobian algebras for tensor algebras of bimodules}}),   we investigate cluster tilting within the framework of $2$-Calabi-Yau or exact stably
$2$-Calabi-Yau categories $\C$ over  an arbitrary field $\kk$.

Simply-laced preprojective algebras are proven useful to get examples of exact stably $2$-Calabi-Yau categories having  cluster tilting subcategories without  loops or $2$-cycles.
A step forward, motivated by this work, could be  a further study of  the representation theory of non simply-laced preprojective algebras \cite{DR80}. A recent account of the study of preprojective algebras over a non algebraically closed field  with respect to cluster tilting is  given in \cite[\S~4.3]{IyamaOppermann2011}.
 
\end{frame}

\section{Some preliminaries on triangulated and exact categories}
\begin{frame}{Framework and conventions}
\begin{itemize}
\only<1>{
\item Fix a ground field $\kk$ and write $\Dual=\HM{\kk}{\sdash,\kk}$ for the  standard duality.
\item An \emph{exact category} is  an additive category   endowed with a set of  \emph{short exact sequences}  denoted by $\Seq{L\to[/above={h}] M\to[/above={f}] N}$ and satisfying the Quillen's axioms, see  \cite{Keller96} and \cite{Buhler}.    An exact category is Frobenius if it has enough projectives and enough injective and the projectives coincide with the injectives.
\item For a triangulated  category $\C$,  we shall often denote a triangle $\trgl{X}{U}{Y}$ simply by $\Seq{X \to U \to Y}$ and when this happens we also omit to specify the fourth component $a[1]$ in a morphism of triangle $(a,u,b,a[1])$. We also   write $\Ext{\C}(X,Y)=\C(X,Y[1])\cong \C(X[-1],Y)$ for all $X,Y\in \C$.
\item Calabi-Yau condition.  Recall that  $\C$ is  \emph{$2$-Calabi-Yau} (weakly $2$-Calabi-Yau in \cite{Keller08}) if for all objects $X,Y \in \C$ we have a functorial isomorphism
$\morph[/isomark={\sim}]{\Dual\Ext{\C}(X,Y)}{\Ext{\C}(Y,X)}$, or equivalently, we have a functorial isomorphism   $\morph[/isomark={\sim}]{\Dual\Hc{X,Y}}{\Hc{Y,X[2]}}$. Recall from \cite{BIRSc} that a  \emph{stably $2$-Calabi-Yau}  category  is a Hom-finite Frobenius  category $\C$ whose stable category $\underline{\C}$, which is a triangulated category according to \cite{H1}, is  $2$-Calabi-Yau.}
\only<2>{
\item \textbf{Assumption}: $\C$   stands for  a Hom-finite Krull-Schmidt $\kk$-category which moreover is assumed to be either  triangulated $2$-Calabi-Yau or exact stably $2$-Calabi-Yau. For all $X\in \C$ indecomposable,  $\kk_X$  is the \emph{residue division algebra} of $\C(X,X)$, in other terms
\begin{cent}
$\kk_X= 
{\C(X,X)}/{\J{\C}(X,X)}$
\end{cent} 
where $\J{\C}$ is the (Jacobson)  \emph{radical bifunctor} of $\C$.
\item By a subcategory of $\C$  we always mean an additive full subcategory of $\C$ which is  stable by direct summands, direct sums and isomorphisms. }
\end{itemize}
 \end{frame}

\begin{frame}{Octaedral axiom  in action in exact categories}

\begin{Lem} \label{lem:exactcat+octaedral}
For any  diagram $(\Delta)$,  $(\Delta')$, $(\op{\Delta})$ or  $(\op{\Delta'})$ below 
\begin{cent}
$\myeqar[cccc]{
(\Delta): & (\Delta'): & (\op{\Delta}): & (\op{\Delta'}):  \\ 
 \diagram[/dist=1,/clsep=1.5]{{M} \& {Z} \& {X} \\
{M} \& {B} \& {N} \\ {} \& {Y} \& {Y} }{
\path[->]
(m-1-1) edge[dashed] node[above] {$h'$} (m-1-2)
(m-1-2) edge[dashed] node[above] {$f'$} (m-1-3)
(m-2-1) edge node[below] {$h $} (m-2-2)
(m-2-2) edge node[below] {$f$} (m-2-3)
(m-1-2) edge[dashed] node[left] {$u'$} (m-2-2)
(m-2-2) edge[dashed]    (m-3-2)
(m-1-3) edge node[right] {$u$} (m-2-3)
(m-2-3) edge node[right] {$v$} (m-3-3);
\path (m-1-1) edge[white] node[black,sloped] {\eqto} (m-2-1);
\path (m-3-2) edge[white] node[black] {\eqto} (m-3-3);
\path (m-1-2) edge[white] node[black] {{\textsc{pb}}} (m-2-3);
}, &
 \diagram[/dist=1,/clsep=1.5]{{} \& {Y} \& {Y} \\ {M} \& {Z} \& {X} \\
{M} \& {B} \& {N}  }{
\path[->]
(m-2-1) edge[dashed] node[above=-0.08] {$h'\,$} (m-2-2)
(m-2-2) edge[dashed] node[above=-0.08] {$f'$} (m-2-3)
(m-3-1) edge node[below] {$h $} (m-3-2)
(m-3-2) edge node[below] {$f$} (m-3-3)
(m-1-2) edge[dashed] node[left=-0.08] {$s'$}    (m-2-2)
(m-1-3) edge node[right] {$s$} (m-2-3)
(m-2-2) edge[dashed] node[left=-0.08] {$t'$} (m-3-2)
(m-2-3) edge node[right] {$t$} (m-3-3);
\path (m-1-2) edge[white] node[black] {\eqto} (m-1-3);
\path (m-2-1) edge[white,sloped] node[black] {\eqto} (m-3-1);
\path (m-2-2) edge[white] node[black] {{\textsc{pb}}} (m-3-3);
},    & 
\diagram[/dist=1,/clsep=1.5]{ {Y} \& {Y} \& {} \\ {N} \& {B} \& {M} \\
{X} \& {Z} \& {M}  }{
\path[->]
(m-2-1) edge node[above] {$\op{f}$} (m-2-2)
(m-2-2) edge  node[above] {$\op{h}$} (m-2-3)
(m-3-1) edge[dashed]   (m-3-2)
(m-3-2) edge[dashed]    (m-3-3)
(m-1-1) edge node[left] {$\op{v}$} (m-2-1)
(m-2-1) edge node[left] {$\op{u}$} (m-3-1)
(m-1-2) edge[dashed]   (m-2-2)
(m-2-2) edge[dashed]  (m-3-2);
\path (m-1-1) edge[white] node[black] {\eqto} (m-1-2);
\path (m-2-3) edge[white,sloped] node[black] {\eqto} (m-3-3);
\path (m-2-1) edge[white] node[black] {{\textsc{po}}} (m-3-2);
}, &  
\diagram[/dist=1,/clsep=1.5]{{N} \& {B} \& {M} \\
{X} \& {Z} \& {M} \\ {Y} \& {Y} \& {} }{
\path[->]
(m-1-1) edge node[above] {$\op{f} $} (m-1-2)
(m-1-2) edge node[above] {$\op{h}$} (m-1-3)
(m-2-1) edge[dashed]  (m-2-2)
(m-2-2) edge[dashed]   (m-2-3)
(m-1-1) edge node[left] {$\op{t}$} (m-2-1)
(m-2-1) edge node[left] {$\op{s}$} (m-3-1)
(m-1-2) edge[dashed]  (m-2-2)
(m-2-2) edge[dashed]    (m-3-2);
\path (m-1-3) edge[white] node[black,sloped] {\eqto} (m-2-3);
\path (m-3-1) edge[white] node[black] {\eqto} (m-3-2);
\path (m-1-1) edge[white] node[black] {{\textsc{po}}} (m-2-2);
}, } $
\end{cent}
 in which  a short exact plain row $\Seq{*\to * \to *}$ and a short exact plain column $\Seq{*\to * \to *}$ are given, the  dashed morphisms exist and complete  the diagram to a commutative diagram in which the dashed row $\Seq[/Dir={->}]{ * \to[style=dashed] * \to[style=dashed] * }$ and the dashed column $\Seq[Dir={->}]{*\to[style=dashed] * \to[style=dashed] *}$   are also short exact. The square marked by \textsc{pb} (resp., \textsc{po}) is a pull-back (resp., push-out) square.
\end{Lem}

\begin{itemize}
\item<2> Notice that $(\op{\Delta})$ expresses  "Noether isomorphism $B/N\cong (B/Y)/(N/Y)$" (\cite[Lem~3.5]{Buhler}).
\end{itemize}
\end{frame}

\begin{frame}{Triangulated structure of  stable categories}
\only<1-2>{ In the exact case ( ou category  $\C$ is then assumed Frobenius), the triangulated structure of $\undl{\C}$ is given as follows. 
\begin{itemize}
\item For all $u\in\C(X,Y)$  we denote by $\bbar{u} \in \undl{\C}(X,Y)$  the residue class of  $u$ with respect to the ideal of morphisms in $\C$  factoring through projective-injective.
\item  Each object $X\in\C$  is associated with a chosen short exact sequence denoted by $\Seq{(\inj[X],\proj[X]): X\to[/above={\inj[X]}]  \, \iinj X \to[/above={\proj[X]}] X[1]}$  with $\iinj X$ a projective-injective object.
\item Each short exact sequence $\Seq{(\xi): X \to[/above={h}] B  \to[/above={f}] Y }$  defines a \emph{standard triangle} $\Seq{(\undl{\xi}) : X \to[/above={\bbar{h}}] B  \to[/above={\bbar{f}}] Y \to[/above={\bbar{\delta}}] X[1] }$  
\item<2>  and for all $a\in \C(X,Z)$,  the shift $\bbar{a}[1]\in\undl{\C}(X[1],Z[1])$  is given by $\bbar{a\prim}$, where  $\delta\in \C(Y,X[1])$   and  $a'\in\C(X[1],Z[1])$  are any maps occurring  as in commutative diagrams of the forms:
\begin{equation} \label{eq:trigstruc.stablecat}
\diagram{ {(\xi):} \&[-0.5em] {X} \& {B} \& {Y} \\ 
{(\inj[X],\proj[X]):} \&[-0.5em] {X} \& {\iinj X} \&  {X[1]} }{
\path[->] (m-1-2) edge node[above]{$h$} (m-1-3) 
(m-1-3) edge node[above]{$f$} (m-1-4)
(m-2-2) edge node[below]{$\inj[X]$} (m-2-3)
(m-2-3) edge node[below]{$\proj[X]$} (m-2-4)
(m-1-3) edge node[left]{$h'$} (m-2-3)
(m-1-4) edge node[right]{$\delta$} (m-2-4);
\path (m-1-2) edge[white] node[black,sloped] {\eqto} (m-2-2); }
\quad \text{ and }
\diagram{ {(\inj[X],\proj[X]):} \&[-0.5em]  {X} 
\& {\iinj X} \&  {X[1].} \\
{(\inj[Z],\proj[Z]):} \&[-0.5em]  {Z} 
\& {\iinj Z} \&  {Z[1].} }{
\path[->] (m-1-2) edge node[above]{$\inj[X]$} (m-1-3) 
(m-1-3) edge node[above]{$\proj[X]$} (m-1-4)
(m-2-2) edge node[below]{$\inj[Z]$} (m-2-3)
(m-2-3) edge node[below]{$\proj[Z]$} (m-2-4)
(m-1-2) edge node[left]{$a$} (m-2-2)
(m-1-3) edge node[left]{$\alpha$} (m-2-3)
(m-1-4) edge node[right]{$a'$} (m-2-4);}
\end{equation}
\end{itemize}
 }
\only<3>{ 
\begin{itemize}
\item We then point out the following observation.
\end{itemize}
\begin{Rem} \label{rem.Fxt:stablecat} Keeping above notations, for all $X,Y\in \C$ we have a canonical isomorphism 
\begin{cent}
$\morph[/isomark={\sim}]{\Ext{\C}(Y,X)}{ \Ext{\undl{\C}}(Y,X)\, : (\xi)\mapsto \bbar{\delta}}$.
\end{cent}
Moreover,  if  $X$   and $Y$   have no projective direct summand, then $\J{\undl{\C}}(X,Y)=\undl{\J{\C}(X,Y)}=\set{\bbar{u}~:~ u\in \J{\C}(X,Y)}$, and in particular, we have a natural identification of residue $\kk$-algebras:
\begin{cent}
$\morph[/isomark={\sim}]{\kk_X}{ \, {\undl{\C}(X,X)}/{\J{\undl{\C}}(X,X)}\, : x +\J{\C}(X,X) \mapsto \bbar{x}+\J{\undl{\C}}(X,X) }$.
\end{cent}  
\end{Rem}
}
\only<4>{
\begin{itemize}
\item As a consequence of previous remark, the $2$-Calabi-Yau condition on $\undl{\C}$ amounts to the following:  we have  functorial isomorphisms
\begin{cent}
$\morph[/isomark={\sim}]{\Dual\Ext{\C}(X,Y)}{\Ext{\C}(Y,X)}$, with  $X,Y \in \C$.
\end{cent}
\end{itemize}
 \begin{itemize}
 \item The following easy fact is needed.
 \end{itemize} 
 \begin{Lem} \label{lem.morph-stablecat} Keeping previous notations, for any endomorphism $(a,b,c)$ of the short exact sequence $(\xi)$, the residue class $(\bbar{a},\bbar{b},\bbar{c},\bbar{a}[1])$ is an endomorphism of the induced standard triangle $(\undl{\xi})$.
 \end{Lem}
 }
\end{frame}

\section{Minimal approximation sequences}

\subsection{Minimal maps, approximation sequences}

\begin{frame}{Minimal maps}

\begin{Defn} \label{def.minmap}
A morphism  $\morph{f: X}{Y}$  in $\C$  is called \emph{right minimal} if for any direct summand $\morph{U}{0}$ of $f$  we have $U=0$. Dually, $f$ is called left minimal if for any direct summand $\morph{0}{U}$ of $f$  we have $U=0$.
\end{Defn}

The following  and its dual version are valid over arbitrary Krull-Schmidt categories. 
\begin{Lem} \label{lem:char.minmap}
For any morphism $\morph{f: X}{Y}$ in $\C$, we have: $f$   right minimal if and only if every section $s$ such that $fs=0$
is zero, if and only if every morphism $u$ such that $fu = 0$ is a radical map, if and only if every $u\in \C(X,X)$ such that $fu=f$ is an automorphism.
\end{Lem} 

The following   also applies to all  Krull-Schmidt categories or $\Hm$-finite categories.
\begin{Lem}\label{lem.min-map} Let $f\in\C(X,Y)$. Then  there are isomorphisms $\morph[/isomark={\sim}]{\gamma: X_1\oplus X_2}{X}$ and $\morph[/isomark={\sim}]{\lambda:Y}{ Y_1\oplus Y_2}$ such that 
$\morph{f\gamma=\Psmatr{f' & 0}:X_1\oplus X_2 }{Y}$  and $\morph{\lambda f=\Psmatr{f'' \\ 0} : Y}{Y_1\oplus Y_2}$ with $f'$ right minimal and $f''$ left minimal.
\end{Lem}

\end{frame}

\begin{frame}{Approximation sequences}
\only<1>{
\begin{itemize}
\item  Let $\cal{M}$ be a subcategory of $\C$ and $B\in \cal{M}$.  A map $f\in \C(B,X)$ is a \emph{right $\cal{M}$-approximation} if  each object $Z\in\cal{M}$ yields an exact sequence $\Seq{\C(Z,B) \to[/above={f \sdash}]  \C(Z,X)  \to  0}$. In addition, it is called minimal if  $f$ is  right minimal. We also have the dual notion of (minimal) \emph{left $\cal{M}$-approximation}.
\item By an \emph{extension} from $Y$  to $X$ in the exact case (resp., in the triangulated case) we  mean any short exact sequence (resp., any triangle) $\Seq{(\xi): X \to[/above={h}]  B\to[/above={f}] Y}$.
\item An  extension $ \Seq{X \to[/above={f}]  Z\to[/above={g}] Y}$  is  called an \emph{$\cal{M}$-approximation sequence} if $f$ is a left $\cal{M}$-approximation and $g$  a right $\cal{M}$-approximation. In addition, it is called \emph{minimal} if $f$ is left minimal and $g$   right minimal.  Notice that minimal approximation sequences never split.  
\end{itemize}
 }
 
\only<2>{ 
\begin{Rem} Assuming the exact case,  let $\undl{\cal{M}}$ be the subcategory of $\undl{\C}$ induced by $\cal{M}$. Then a short exact sequence  $ \Seq{(\xi) : X \to[/above={f}]  Z\to[/above={g}] Y}$  is an $\cal{M}$-approximation sequence (resp.,  a minimal $\cal{M}$-approximation sequence) if and only if the induced standard triangle $(\undl{\xi})$ is an $\undl{\cal{M}}$-approximation sequence (resp., a minimal $\undl{\cal{M}}$-approximation sequence). 
\end{Rem}

A direct application of previous  Lemma  yields the following result.  

\begin{Lem}\label{lem.min-seq}  Any extension $\Seq{(\xi): X \to[/above={h}] B\to[/above={f}] Y}$ in $\C$   is a direct sum of a minimal extension
 $ \Seq{(\xi_{\mrm{min}}): X_1\to[/above={h_1}]  B_1\to[/above={f_1}]  Y_1}$  and a split  extension  
 $\Seq{ (\xi_0): X_2\to[/above={\Psmatr{1\\ 0}}]  X_2\oplus Y_2\to[/above={\Psmatr{0 & 1}}]  Y_2}$.
\end{Lem}

\begin{Cor}\label{cor.minseq} For any non-split extension $\Seq{(\xi): X \to[/above={h}]  B\to[/above={f}] Y}$  in $\C$, if $Y$ is indecomposable then $h$ is a left minimal radical morphism, if $X$ is indecomposable then $f$ is a right minimal
radical morphism.  
\end{Cor}
}
\end{frame}

\begin{frame}{Minimality and Rigid  subcategories}
\only<1>{ 
We recall (for example from \cite{BMRRT,DK,Iyama2007}) that a subcategory $\T\subset\C$ is called \emph{rigid} if $\Ext{\C}(X,Y)=0$   for all $X,Y\in \T$. It is called \emph{cluster tilting} if it is \emph{functorially finite} (i.e. all $X\in \C$ admits left and right $\T$-approximations) and
\begin{cent}
$\T=\set{X\in\C~:~\forall T\in\T,\Ext{\C}(X,T)=0}$ (equivalently, $\T=\set{X\in\C~:~\forall T\in\T,\Ext{\C}(T,X)=0}$).
\end{cent}
It is immediate that a cluster tilting subcategory is also rigid.}
\only<2>{  
The  following consequence is also derived.

\begin{Lem} \label{lem:approxseq.thesame-end-terms} Suppose that $\cal{M}\subset \C$ is rigid and $\Seq{(\xi): X\to[/above={h}] U\to[/above={f}] Y }$  and $\Seq{(\xi'): Y\to[/above={h'}] U'\to[/above={f'}]  X}$ are $\cal{M}$-approximation sequences.  Then there is a common direct summand $X_2\oplus Y_2$ of $U$  and $U'$  such that $(\xi)$ is the direct sum of a minimal $\cal{M}$-approximation sequence $\Seq{(\xi_{\mrm{min}}): X_1\to[/above={h_1}] U_1\to[/above={f_1}] Y_1 }$   with the  split extension $X_2\to X_2\oplus Y_2 \to Y_2$, and $(\xi')$ is the direct sum of a minimal  $\cal{M}$-approximation sequence  $\Seq{(\xi_{\mrm{min}}'): Y_1\to[/above={h_1'}] U_1'\to[/above={f_1'}]  X_1}$  with the split extension 
$Y_2\to Y_2\oplus  X_2 \to X_2$.
\end{Lem}
}

\end{frame}

\begin{frame}{From minimal approximation sequences to isomorphisms of Division algebras}
We also prove the following lemma, crucial for the rest of this section.
\begin{Lem} \label{lem.min-approxseq} Let  $\cal{M} $ be a subcategory of $\C$.  If $\Seq{(\xi): X \to[/above={h}] B \to[/above={f}] Y }$ and 
$\Seq{(\xi'): X' \to[/above={h'}] B'  \to[/above={f'}] Y' }$ are minimal $\cal{M}$-approximation sequences then the following assertions hold.
\begin{enumerate}
\item[$\bullet$] Any morphism $\morph{u: X }{X'}$ (resp., $\morph{v: Y }{Y'}$) extends to a morphism $\morph{\varphi: (\xi) }{(\xi')}$ of  approximation sequences. Moreover, $u$ (resp., $v$) is an isomorphism if and only if $\varphi$ is an
isomorphism. 
\item[$\bullet$] $X$ is indecomposable if and only if $Y$ is. In this case,  $(\xi)$ induces an isomorphism $\morph[/isomark={\sim}]{\phiso[\xi]: \kk_X }{\kk_{Y}}$ taking $u+\J{\C}(X,X)$ to $v+\J{\C}(Y,Y)$ whenever the pair $(u,v)$ extends to an endomorphism of $(\xi)$.  
\end{enumerate}
\end{Lem}

Specializing to the exact setting, we also prove the following fact.
\begin{Lem} \label{lem:approxseq.approxstdtriang+induceiso} Assuming the exact case, let  $\cal{M} $ be a subcategory of $\C$ and  $\Seq{(\xi): X \to[/above={h}] B \to[/above={f}] Y }$ a minimal $\cal{M}$-approximation sequence with $X$  and $Y$ indecomposable non-projective. Then for the induced minimal  $\undl{\cal{M}}$-approximation standard triangle $(\undl{\xi})$ we have   $\phiso[\xi]=\phiso[\undl{\xi}]$.
\end{Lem}
\end{frame}

\subsection{Exchange sequences and irreducible maps}

\begin{frame}{Irreducible maps}
\begin{itemize}
\only<1>{ \item We let $\T$ be a subcategory of $\C$ and remind the reader that, by assumption, $\T$ is full, stable by direct summands, direct sums and isomorphisms. Given any subcategory $\cal{X}\subset\T$,  the \emph{factor category} $\T/\cal{X}$  is  the subcategory of $\T$   whose non-zero objects  consist of    all objects   of $\T$   not belonging  to $\cal{X}$. In case $\cal{X}=\add(X)$  for some $X\in \T$,  we   simply  write  $\T/X$  for $\T/\cal{X}$.
Let  $X,Y\in  \T$.   We have $\J{\T}(X,Y)=\J{\C}(X,Y)$, and the square $\J[2]{\T}$ of the radical bifunctor of $\T$ is such that
\begin{cent}
$\J[2]{\T}(X,Y)=\Som{Z\in \T}\J{\C}(Z,Y)\J{\C}(X,Z)=\set{vu~:~\exists Z\in \T, u\in\J{\C}(X,Z),v\in \J{\C}(Z,Y) }$.
\end{cent}
\item Next, define the \emph{$\kk_Y\sdash\kk_X$-bimodule of $\T$-irreducible morphisms} from $X$ to $Y$  to be
\begin{cent}
$\Irr[\T](X,Y)={\J{\C}(X,Y)}/{\J[2]{\T}(X,Y)}$.
\end{cent}
For an indecomposable $X\in \T$, the category $\T$ is said to have  \emph{no loop at $X$} if $\Irr[\T](X,X)=0$.  
}
\only<2>{
\item 
We therefore get the following. 
\begin{Rem} \label{rem.noloop}
 The category $\T$ has no loop at the indecomposable object  $X$ if and only if any radical morphism $u\in \J{\C}(X,X)$ factors through some object in $\T/X$.
 Consequently, in this case,  any left (resp., right)  $(\T/X)$-approximation $\morph{h:X}{B}$ (resp., $\morph{f: B}{X}$) is left almost split (resp., right almost split).
 \end{Rem}
}
\end{itemize}
\end{frame}

\begin{frame}{Exchange sequences}
\begin{Defn} \label{def:exchseq}
For a subcategory $\cal{M}$ of $\C$, a minimal $\cal{M}$-approximation sequence $\Seq{(\xi):X\to B \to Y}$ is called an \emph{exchange sequence} if $X$ is indecomposable, $X\not\in\cal{M}$ and $\add(X, \cal{M})$ is rigid, or, equivalently, $Y$ is indecomposable, $Y\not\in\cal{M}$ and $\add(Y, \cal{M})$ is rigid. More precisely, $(\xi)$ is called an exchange sequence from  $\add(X, \cal{M})$  to $\add(Y, \cal{M})$.
\end{Defn}

\begin{Lem}\label{lem.exch-seq}  Let $\Seq{(\xi): X  \to[/above={h}]  B  \to[/above={f}]  Y}$  be an exchange sequence from a rigid subcategory $\T$ to a rigid subcategory $\T'$
of $\C$. Then the following are equivalent:
\begin{enumerate}
\item[$\msf{(a)}$] The category $\T$ has no loop at $X$.
\item[$\msf{(b)}$] The category $\T'$ has no loop at $Y$.
\item[$\msf{(c)}$] $\Ext{\C}(Y,X)\cong \kk_X\cong \kk_Y$. 
\end{enumerate}
 Moreover, in this case, any non-split extension from $Y$ to $X$ is isomorphic to  $(\xi)$.
\end{Lem}

\end{frame}

\begin{frame}{Trace maps play a determinant role}
\only<1>{
A $\kk$-linear \emph{trace} on $\K$ is a central element  of the standard dual bimodule $\Dual(\K)$.
\begin{Rem} \label{rem:trace+dualisingbimod} Any finite-dimensional division $\kk$-algebra has a non-zero trace. In particular,  for finite-dimensional division $\kk$-algebras $\K$  and $\K'$  and for any $\K\sdash\K'$-bimodule $B$ (finite-dimensional over $\kk$), the choice of two non-zero traces $\tr\in \Dual(\K)$  and $\tr'\in\Dual(\K')$  induces canonical isomorphisms  
\begin{cent}
$\morph[/above={\sim}]{\Hm[\K](B,\K)}{\Dual(B):\, u\mapsto \tr\circ u}$, ~~~~ 
$\morph[/above={\sim}]{\Hm[\K'](B,\K')}{\Dual(B):\, u\mapsto \tr'\circ u}$
\end{cent} 
 from the left dual (resp., right dual) of $B$ to its standard dual.
\end{Rem}

We now prove the following  connection between   exchange sequences and bimodules of irreducible maps.
\begin{Prop}\label{prop:exchseq+Irred} Let  $\Seq{(\xi): M   \to[/above={h}]  B \to[/above={f}] M'}$ be    an exchange sequence between rigid subcategories $\T$  and $\T'$ and assume that  $\T$ has no loop at $M$. Under the canonical isomorphism $\morph{\phiso[\xi]: \kk_M}{\kk_{M'}}$, for  all $X\in\T\cap\T'$ we have an isomorphism of $\kk_M\sdash\kk_X$-bimodules 
\begin{cent}
$\morph[/above={\sim}]{ \phiso[\xi,X]: \Irr[\T'](X,M')}{\Dual\Irr[\T](M,X)}$.
\end{cent}
\end{Prop}
}
\only<2-4>{
\begin{prv}
\begin{itemize}
\only<2>{ 
\item  Let $X\in\cal{T}\cap\cal{T}'$ be  indecomposable, put $\Bimd{M}{X}=\Irr[\T](M,X)$ and $\Bimd{X}{M'}=\inddl{\Irr[\T'](X,M')}{\phi}{} $. For all $u\in\J{\C}(M,X)$  and $v\in \J{\C}(X,M')$, put $\wtilde{u}=u+\J[2]{\T}(M,X)$  and $\wtilde{v}=v+\J[2]{\T'}(X,M')$.
}
\only<3>{
\item 
 Thus, in view of Remark~\ref{rem.noloop},  $h$ is left minimal almost split in $\T$  while $f$ is right minimal almost split in $\T'$. Hence, for  any two radical maps  $u\in\J{\C}(M,X)$  and $v\in \J{\C}(X,M')$ as before
 and   all $a\in\C(M,M)$, invoking the factorization property of $h$ and $f$ together with Lemma~\ref{lem.min-approxseq}  we get a commutative diagram with the following shape:
\begin{equation}\tag{D1} \label{eq:diag.prv}
\diagram[/dist=1,/clsep=3.5]{ \&[-2em] \&  {X} \&   \\ 
(\xi): \&[-2em] {M} \& {B} \& {M'} \\ 
(\xi): \&[-2em] {M} \& {B} \& {M'} \\ 
\&[-2em] \& {X} \& }{
\path[->]
(m-2-2) edge node[above]  {$h$} (m-2-3)
(m-2-3) edge node[above]  {$f$} (m-2-4)
(m-3-2) edge node[above]  {$h$} (m-3-3)
(m-3-3) edge node[above]  {$f$} (m-3-4)
(m-3-2) edge node[below] {$u$} (m-4-3)
(m-1-3) edge node[above] {$v$}  (m-2-4)
(m-2-2) edge node[left] {$a$}  (m-3-2)
(m-2-4) edge[dashed] node[right] {$a'$}  (m-3-4)
(m-1-3) edge[dashed] node[left] {$v_*$} (m-2-3)
(m-2-3) edge[dashed] node[right] {$b$} (m-3-3)
(m-3-3) edge[dashed] node[right] {$u_*$} (m-4-3);
}
\end{equation}
}
\only<4>{
\item  We claim that there is an isomorphism of $\kk_M\sdash\kk_X$-bimodules given by 
\begin{cent}
$\vphi: \Bimd{X}{M'} \, \to \, \Hm[\kk_X](\Bimd{M}{X},\kk_X) :~
\wtilde{v} \mapsto (\wtilde{u} \mapsto \vphi(\wtilde{v})(\wtilde{u}) =\bbar{u_*v_*})$.
\end{cent} 
................................
} 
\end{itemize}
\end{prv}
}
\end{frame}


\section{Main statements}
\begin{frame}{Comaptible exchange sequences and generalized cluster structures}
\only<1>{
To each subcategory $\T$ of $\C$ can be associated the following combinatorial data.

\begin{Defn} The \emph{exchange matrix} of $\T$ is the (not necessarily finite) matrix $\B(\T)$ which
has coefficients 
\begin{cent}
$\bss[X\nsc Y]=\dim_{\kk_X}\Irr[\T](X,Y)-\dim_{\kk_X}\Irr[\T](Y,X)$
\end{cent} for all $X,Y\in \Ind\T$, $Y$  being non-projective (in the exact case).
\end{Defn}

  The principal part of the exchange matrix $\B(\T)$  is skew symmetrizable via the
diagonal matrix with coefficients $\dim_{\kk}(\kk_X)$ for $X$ indecomposable non-projective. 

Notice that the exchange matrix $\B(\T)$ uniquely determines the valued quiver (or the quiver in the simply-laced case) of $\T$.}

\only<2>{  Now let $M\in\Ind\T$ and suppose that  we have two exchange sequences  $\Seq{(\xi): M^* \to B \to M}$ and $\Seq{(\xi'): M \to B' \to M^*}$  relating  $\T$ to a subcategory $\T^*= \add(\T/M, M^*)$ for some indecomposable object  $M^*\not\in\T$. When $\T$   and $\T^*$  are cluster tilting, it is proved in the simply-laced context of \cite[Theorem~1.6]{BIRSc} that $\B(\T)$ is related to $\B(\T^*)$ via the Fomin-Zelevinsky mutation \cite{FZ1}. But  when the base field $\kk$ is not algebraically closed, the matrix $\B(\T)$ alone does not give enough information about irreducible maps in $\T$, the  corresponding information is  encoded by the modulated quiver of $\T$, (see for example \cite{DR} for the notion of modulated quiver). 

Then,  we introduce the following  notion of compatibility in order to get a more finer relation between  bimodules of irreducible maps in $\T$  and in $\T^*$.   

\begin{Defn} \label{def:compatible.exchseq} The exchange sequences $(\xi)$  and $(\xi')$  are called \emph{compatible} if the induced isomorphisms $\Seq{\phiso[\xi]: \kk_{M^*} \to \kk_M }$ and $\Seq{ \phiso[\xi']: \kk_M \to \kk_{M^*}}$  and   are inverse from each other.
\end{Defn}
}
\only<3>{ 
The following proposition sheds more light on the compatibility between exchange sequences $(\xi)$  and $(\xi')$; they may fail to be compatible especially when the residue division algebra $\kk_M$ is not commutative.

\begin{Prop} \label{prop.compatibility-exch-seq}
 Suppose that $\Seq{(\eta): M \to B'' \to M^*}$ is another exchange sequence from $\T$  to $\T^*$.
Then any $a \in \C(M, M)$ extends to endomorphisms $\morph[/dist=1]{(c,b,a):(\xi)}{(\xi)}$, $\morph[/dist=1]{(a,b',c'):(\xi')}{(\xi')}$ and $\morph[/dist=1]{(a,u^{-1}b'u,\lambda^{-1}c'\lambda):(\eta)}{(\eta)}$ such that $c-c'\in\J{\C}(M^*,M^*)$,
$\lambda$ is an automorphism of $M^*$ and $u\in\C(B',B'')$
is an isomorphism. Thus, when $\kk_M$ is commutative, any two exchange sequences  $\Seq{ M^* \to U \to M}$  and $\Seq{ M \to U' \to M^*}$ relating $\T$ and $\T^*$  are compatible. This is the case if  $\kk$ is algebraically closed.
\end{Prop}

\ParIt{Cluster structure} 
The notion of cluster structure in the present framework  is a natural generalization of the one from \cite{BIRSc}, 
 along the above lines.  Mutation of quivers from \cite{BIRSc} should be replaced by mutation of skew symmetrizable matrices  and exchange sequences are required to be compatible. Recall  that in the simply laced context of \cite{BIRSc}, exchange sequences are automatically compatible according to previous  Proposition.
}
\end{frame}

\begin{frame}{Main statements}
\only<1>{ 
\begin{Thm}[{\cite[Theorem 5.3]{IyamaYoshino08} Iyama-Yoshino}] \label{theo.exch-seq:IY} Let $\cal{M}$ be a subcategory of $\C$ and $M\in \C \ssminus  \cal{M} $ be indecomposable (non-projective in the exact case), such that $\T = \add(X, \cal{M})$ is cluster tilting. Then there is  a unique (up to isomorphism)  $M^*\in \C \ssminus \T$ indecomposable such that 
the subcategory $\mu_M(\T)= \add(M^*,  \cal{M})$ is also cluster
tilting. Moreover, there exist exchange sequences $\Seq{ (\xi): ~ M^* \to B \to M}$ and $\Seq{ (\xi'):  ~ M \to B' \to M^*}$  relating $\T$  and $\mu_M(\T)$.
\end{Thm}

The subcategory $\T^*=\mu_M(\T)$ is called the mutation of $\T$ at $M$ and we also have $\T=\mu_{M^*}(\T^*)$.

Our first contribution then specializes Iyama-Yoshino Theorem  in this context  as follows.

\begin{Prop} \label{prop.compatible-exch-seq}
If $\T$ has no loop at $M$, then  exchange sequences $(\xi)$ and $(\xi')$ can be chosen to be compatible.
\end{Prop} 
}
\only<2>{ 

Before stating the next result of this paper,  recall that $\T$   has
no $2$-cycle if for any $X,Y\in \T$ indecomposable, either $\Irr[\T](X,Y)=0$, either $\Irr[\T](Y,X)=0$. 
\begin{Thm}\label{theo.clstr}
 If $\T$ is a cluster tilting subcategory of $\C$ without loop or  $2$-cycle, and   $\T^*$  the mutation of $\T$  at some  
 $M\in\Ind\T$ 
 then \vskip-1.5em
 \begin{cent}
 $\B(\T^*)=\mu_M(\B(\T))$.
 \end{cent}  
\end{Thm}

We should emphasize that the main contribution of this work is two fold: 
\begin{itemize}
\item Our main  Theorem  compares exchange matrices (or equivalently, valued quivers) of $\T$  and $\T^*$ using Fomin-Zelevinsky mutation.
\item The existence of compatible exchange sequences given  by the Proposition more above  is crucial fro describing  bimodules of irreducible maps in $\T^*$,  indeed this result  yields a first step  to describe the modulated quiver of   $\T^*$ using  that of $\T$.\end{itemize}     
 }
 \only<3>{
 \begin{itemize}
 \item We get the following consequence where we use 
  \cite[Proposition~6.14]{BMRRT} for the second part.
 \end{itemize} 
     
  \begin{Cor} \label{cor.clustcat:clst}
   Suppose that none of cluster tilting subcategories of $\C$ have loop or $2$-cycle. Then $\C$ has a cluster structure induced by its cluster tilting subcategories.
  In particular, it is the case for cluster categories associated with  finite dimensional hereditary algebras. 
  \end{Cor}
 }
\end{frame}

\bibliographystyle{elsarticle-num}

\end{document}